\documentclass{article}
\usepackage[usenames]{color}
\usepackage{amssymb}
\usepackage{amsmath, amsthm}
\usepackage{amscd}
\usepackage{fancyhdr}
\usepackage[latin1]{inputenc}
\usepackage[T1]{fontenc}

\def\R{\mathbb{R}}

\def\E{\mathbb{E}}

\newcommand{\F}{{\cal F}}

\newcommand{\X}{{\cal X}}
\newcommand{\D}{{\cal D}}

\newcommand{\1}{{\rm 1}\kern-0.24em{\rm I}}

\def\A{{\cal A}}

\def\F{{\cal F}}

\def\Z{{\cal Z}}

\def\Ec{{\cal E}}

\def\Y{{\cal Y}}
\def\Prob{\mathbb{P}}

\title{High-dimensional stochastic optimization with the generalized Dantzig estimator}

\author{Karim Lounici}

\newtheorem{lem}{Lemma}

\newtheorem{theo}{Theorem}

\newtheorem{cor}{Corollary}
\newtheorem{assum}{Assumption}

\begin{document}
\maketitle

\begin{abstract}
We propose a generalized version of the Dantzig selector. We show
that it satisfies sparsity oracle inequalities in prediction and
estimation. We consider then the particular case of
high-dimensional linear regression model selection with the Huber
loss function. In this case we derive the sup-norm convergence
rate and the sign concentration property of the Dantzig
estimators under a mutual coherence assumption on the dictionary.
\end{abstract}



\section{Introduction}

Let $\Z=\X\times \Y$ be a measurable space. We observe a set of
$n$ i.i.d. random pairs $Z_{i}=(X_{i},Y_{i})$, $i=1,\ldots,n$
where $X_{i}\in\X$ and $Y_{i}\in\Y$. Denote by $P$ the joint
distribution of $(X_{i},Y_{i})$ on $\X\times \Y$, and by $P^{X}$
the marginal distribution of $X_{i}$. Let $Z=(X,Y)$ be a random
pair in $\Z$ distributed according to $P$. For any real-valued
function $g$ on $\X$, define
$||g||_{\infty}=\mathrm{ess}\sup_{x\in\X}|g(x)|$, $\|g\|=\left(
\int_{\X}g(x)^{2}P^{X}(dx) \right)^{1/2}$ and $||g||_{n}  =
\left(  \frac{1}{n}\sum_{i=1}^{n} g(X_{i})^{2} \right)^{1/2}$. Let
$\D=\{f_{1},\ldots,f_{M}\}$ be a set of real-valued functions on
$\X$ called the dictionary where $M\geqslant 2$. We assume that the functions of the dictionary are normalized, so that $\|f_{j}\|=1$ for all $j=1,\ldots,M$. We also assume that $||f_{j}||_{\infty}\leqslant L$ for some $L>0$. For any $\theta
\in \R^{M}$, define $f_{\theta}=\sum_{j=1}^{M}\theta_{j}f_{j}$
and $J(\theta)=\{j\; :\; \theta_{j}\neq 0\}$. Let
$M(\theta)=|J(\theta)|$ be the cardinality of $J(\theta)$ and
$\vec{\mathrm{sign}}(\theta)=(\mathrm{sign}(\theta_{1}),\ldots,\mathrm{sign}(\theta_{M}))^{T}$
where
\begin{equation*}
\mathrm{sign}(t)=\begin{cases}
1 &\text{if $t>0$},\\
0 &\text{if $t=0$},\\
-1 &\text{if $t<0$}.
\end{cases}
\end{equation*}
For any vector $\theta\in\R^{M}$ and any subset $J$ of
$\{1,\ldots,M\}$, we denote by $\theta_{J}$ the vector in $\R^{M}$
which has the same coordinates as $\theta$ on $J$ and zero
coordinates on the complement $J^{c}$ of $J$. For any integers $1
\leqslant d,p < \infty$ and $w=(w_{1},\ldots,w_{d})\in \R^{d}$, the
$l_{p}$ norm of the vector $w$ is denoted by
$|w|_{p}\stackrel{\Delta}{=}\left( \sum_{j=1}^{d}|w_{j}|^{p}
\right)^{1/p}$, and
$|w|_{\infty}\stackrel{\Delta}{=}\max_{1\leqslant j \leqslant
d}|w_{j}|$.

Consider a function $\gamma:\,\Y\times \R \rightarrow \R^{+}$
such that for any $y$ in $\Y$ and $u,u'$ in $\R$ we have
\begin{equation*}
|\gamma(y,u)-\gamma(y,u')| \leqslant |u-u'|.
\end{equation*}
We assume furthermore that $\gamma(y,\cdot)$ is
convex and differentiable for any $y\in\Y$. We assume that for any $y\in\Y$ the
derivative $\partial_{u}\gamma(y,\cdot)$ is absolutely
continuous. Then $\partial_{u}\gamma(y,\cdot)$ admits a
derivative almost everywhere which we denote by
$\partial^{2}_{u}\gamma(y,\cdot)$. Consider the loss function
$Q:\,\Z\times \R^{M} \rightarrow \R^{+}$ defined by
\begin{equation}
Q(z,\theta) = \gamma(y,f_{\theta}(x)).
\end{equation}

The expected and empirical risk measures at point $\theta$ in
$\R^{M}$ are defined respectively by
\begin{equation*}
R(\theta)\stackrel{\triangle}{=}\E\left(Q(Z,\theta)\right),
\end{equation*}
where $\E$ is the expectation sign, and
\begin{equation*}
\hat{R}_{n}(\theta)\stackrel{\triangle}{=}\frac{1}{n}\sum_{i=1}^{n}Q(Z_{i},\theta).
\end{equation*}

Define the target vector as a minimizer of $R(\cdot)$ over $\R^{M}$:
\begin{equation*}
\theta^{*}\stackrel{\triangle}{=}\arg\min_{\theta\in\R^{M}}R(\theta).
\end{equation*}
Note that the target vector is not necessarily unique. From now
on, we assume that there exists a $s$-sparse solution $\theta^{*}$, i.e., a solution
with $M(\theta^{*})\leqslant s$, and that this sparse solution is
unique. We will see that this is indeed the case under the coherence
condition on the dictionary (cf. Section 3 below).

Define the excess risk of the vector $\theta$ by
\begin{equation*}
\Ec(\theta)=R(\theta)-R(\theta^{*}),
\end{equation*}
and its empirical version by
$$\Ec_{n}(\theta)=R_{n}(\theta)-R_{n}(\theta^{*}).$$
Our goal is to derive sparsity oracle inequalities for the excess
risk and for the risk of $\theta^{*}$ in the $l_{1}$ norm and in
the sup-norm.

We consider the following minimization problem:
\begin{align}\label{Dantso}
\min_{\theta\in \Theta}  |\theta|_{1}\hspace{0.2cm} &\text{subject
to}\hspace{0.2cm}
\left|\nabla
\hat{R}_{n}(\theta)\right|_{\infty}\leqslant
r,
\end{align}
where $\nabla\hat{R}_{n}\stackrel{\triangle}{=}
(\partial_{\theta_{1}}\hat{R}_{n},\ldots,\partial_{\theta_{M}}\hat{R}_{n})^{T}$,
$r>0$ is a tuning parameter defined later and $\Theta$ is a
convex subset of $\R^{M}$ specified later. Solutions of
(\ref{Dantso}), if they exist, will be taken as estimators of
$\theta^{*}$. Note that we will prove in Lemma 3 that under
Assumption \ref{A?-norm-sup} the set $ \{\theta \in \Theta \, :
\, \left|\nabla\hat{R}_{n}(\theta)\right|_{\infty}\leqslant r \}$
is non-empty with probability close to one. Note also that in the
applications considered in Section 3, the constraint $|\nabla
\hat{R}_{n}(\theta)|_{\infty}\leqslant r $ can be defined as a
system of inequalities involving convex functions. Thus,
solutions to (\ref{Dantso}) exist and can be efficiently computed
via convex optimization. In particular, for the regression model
with the Huber loss, the gradient $\nabla \hat{R}_{n}(\theta)$ is
piecewise linear so that (\ref{Dantso}) reduces in this case to a
standard linear programming problem. Denote by $\hat{\Theta}$ the
set of all solutions of (\ref{Dantso}). For the reasons above, we
assume from now on that $\hat{\Theta}\neq \emptyset$ with
probability close to one.

The definition of our estimator (\ref{Dantso}) can be motivated as
follows. Since the loss function $Q(z,\cdot)$ is convex and
differentiable for any fixed $z\in\Z$, the expected risk $R$ is
also a convex function of $\theta$ and it is differentiable under
mild conditions. Thus, minimizing $R$ is equivalent to finding the
zeros of $\nabla R$. The quantity $\nabla \hat{R}_{n}(\theta)$ is
the empirical version of $\nabla R(\theta)$. We choose the
constant $r$ such that the vector $\theta^{*}$ satisfies the
constraint $|\nabla \hat{R}_{n}(\theta^{*})|\leqslant r$ with
probability close to $1$. Then among all the vectors satisfying
this constraint, we choose those with minimum $l_{1}$ norm. Note
that if we consider the linear regression problem with the
quadratic loss, we recognize in (\ref{Dantso}) the Dantzig
minimization problem of Candes and Tao \cite{CT07}. From now on,
we will call (\ref{Dantso}) the generalized Dantzig minimization
problem.

Bickel et al. \cite{BRT07}, Candes and Tao \cite{CT07} and
Koltchinskii \cite{Kolt07} proved that the Dantzig estimator
performs well in high-dimensional regression problems with the
quadratic loss. In particular they proved sparsity oracle
inequalities on the excess risk and the estimation of
$\theta^{*}$ for the $l_{p}$ norm with $1\leqslant p \leqslant
2$.

The problem (\ref{Dantso}) is closely related to the minimization
problem:
\begin{equation}\label{lassosol}
\min_{\theta\in\Theta} \hat{R}_{n}(\theta)+r|\theta|_{1},
\end{equation}
which is a generalized version of the Lasso. For the Lasso
estimator, Bunea et al \cite{BTW07b} proved similar results in
high-dimensional regression problems with the quadratic loss
under a mutual coherence assumption \cite{DET06} and Bickel et al
\cite{BRT07} under a weaker Restricted Eigenvalue assumption.
Koltchinskii \cite{Kolt06} derived similar results for the Lasso
in the context of high-dimensional regresssion with twice
differentiable Lipschtiz continuous loss functions under a
restricted isometry assumption. Van de Geer \cite{VdG07,VdG07b}
obtained similar results for the Lasso in the context of
generalized linear models with Lipschtiz continuous loss
functions. Lounici \cite{L08} derived sup-norm convergence rates
and sign consistency of the Lasso and Dantzig estimators in a
high-dimensional linear regression model with the quadratic loss
under a mutual coherence assumption.

The paper is organized as follows. In Section 2 we derive sparsity
oracle inequalities for the excess risk and for estimation of
$\theta^{*}$ for the generalized Dantzig estimators defined by
(\ref{Dantso}) in a stochastic optimization framework. In section
3 we apply the results of Section 2 to the linear regression model
with the Huber loss and to the logistic regression model. In
Section 4 we prove the variable selection consistency with rates
under a mutual coherence assumption for the linear regression
model with the Huber loss. In section 5 we show a sign
concentration property of the thresholded generalized Dantzig
estimators for the linear regression model with the Huber loss.

\section{Sparsity oracle inequalities for prediction and estimation with the $l_{1}$ norm}

We need an assumption on the dictionary to derive prediction and
estimation results for the generalized Dantzig estimators. We first state the Restricted Eigenvalue assumption \cite{BRT07}.
\begin{assum}\label{A1-RE}
\begin{equation*}
\zeta(s) \stackrel{\triangle}{=}
\min_{J_{0}\subset\{1,\ldots,M\}:|J_{0}|\leqslant
s}\hspace{0.5cm}\min_{\Delta\neq
0:|\Delta_{J_{0}^{c}}|_{1}\leqslant
|\Delta_{J_{0}}|_{1}}\frac{||f_{\Delta}||}{|\Delta_{J_{0}}|_{2}}>0.
\end{equation*}
\end{assum}
\noindent It implies an "equivalence" between the two norms
$|\Delta|_{2}$ and $\|f_{\Delta}\|$ on the subset $\{\Delta\neq
0:|\Delta_{J(\Delta)^{c}}|_{1}\leqslant
|\Delta_{J(\Delta)}|_{1}\}$ of $\R^{M}$.

We need the following assumption on $\|f_{\theta^{*}}\|_{\infty}$.
\begin{assum}\label{A?-norm-sup}
There exists a constant $K>0$ such that
$\|f_{\theta^{*}}\|_{\infty}\leqslant K$.
\end{assum}
From now on we take for $\Theta$ the set
\begin{equation*}
\Theta = \{\theta\in\R^{M}\, : \, \|f_{\theta}\|_{\infty}\leqslant
K\}.
\end{equation*}

The following assumption is a version of the margin condition
(cf. \cite{T04}). It links the excess risk to the functional norm
$\|\cdot\|$.

\begin{assum}\label{A2-MA}
For any $\theta\in\Theta$ there exits a constant $c>0$ depending
possibly on $K$ such that
\begin{equation*}
\|f_{\theta}-f_{\theta^{*}}\| \leqslant
c(R(\theta)-R(\theta^{*}))^{1/\kappa},
\end{equation*}
where $1< \kappa \leqslant 2$.
\end{assum}
We will prove in Section 2.1 below that this condition is always satisfied with the constant $\kappa=2$ for the regression model with Huber loss and for the logistic regression model. We also
need the following technical assumption.
\begin{assum}\label{A?-techass}
The constants $K$ and $L$ satisfy
\begin{equation*}
1 \leqslant K,L \leqslant \sqrt{\frac{n}{\log M}}.
\end{equation*}
\end{assum}
Define the quantity
\begin{equation}\label{rtilde}
\tilde{r} = 4\sqrt{2}L\frac{\log M}{n}  + 2 \sqrt{6}\sqrt{\frac{\log
M}{n}}.
\end{equation}
We assume from now on that $\tilde{r}\leqslant 1$.

The main results of this section are the following sparsity oracle
inequalities for the excess risk and for estimation of
$\theta^{*}$ in the $l_{1}$ norm.
Define
\begin{equation}\label{r}
r=6\|\partial_{u}\gamma\|_{\infty}\tilde{r}.
\end{equation}
\begin{theo}
Let Assumptions \ref{A1-RE} - \ref{A?-techass} be satisfied. Take
$r$ as in (\ref{r}). Assume that
$M(\theta^{*})\leqslant s$. Then, with probability at
least $1-M^{-1} -M^{-K} - 3M^{-2K}\log \frac{n}{\log M}$, we have
\begin{equation}\label{th1-pred}
\sup_{\hat{\theta}\in \hat{\Theta}}\Ec(\hat{\theta}) \leqslant
\left(
\frac{2(1+2K)cr\sqrt{s}}{\zeta(s)}\right)^{\frac{\kappa}{\kappa-1}}
+
12\|\partial_{u}\gamma\|_{\infty}\frac{\kappa}{\kappa-1}\tilde{r}^{2},
\end{equation}
and
\begin{equation}\label{th1-est}
\sup_{\hat{\theta}\in \hat{\Theta}}|\hat{\theta}-\theta^{*}|_{1}
\leqslant \left(
\frac{2c\sqrt{s}}{\zeta(s)}\right)^{\frac{\kappa}{\kappa-1}}((1+2K)r)^{\frac{1}{\kappa-1}}
+ \frac{2K}{(\kappa-1)(1+2K)}\tilde{r}.
\end{equation}
\end{theo}

Note that the regularization parameter $r$ does not depend on the
variance of the noise if we consider the regression model with non-quadratic loss. In this case, the use of Lipschtiz losses
enables us to treat cases where the noise variable does not admit a
finite second moment, e.g., the Cauchy distribution. The price to
pay is that we need to assume that
$\|f_{\theta^{*}}\|_{\infty}\leqslant K$ with known $K$.

\begin{proof}
For any $\hat{\theta}\in \hat{\Theta}$ define
$\Delta=\hat{\theta}-\theta^{*}$. We have
\begin{eqnarray}\label{interm1}
\Ec(\hat{\theta}) &\leqslant& \Ec_{n}(\hat{\theta}) +
\Ec(\hat{\theta}) - \Ec_{n}(\hat{\theta})\nonumber\\
&=& \Ec_{n}(\hat{\theta}) + \frac{\Ec(\hat{\theta}) -
\Ec_{n}(\hat{\theta})}{|\Delta|_{1}+\tilde{r}}(|\Delta|_{1}+\tilde{r})\nonumber\\
&\leqslant& \Ec_{n}(\hat{\theta}) + \sup_{\theta \in \Theta :
\theta\neq \theta^{*} } \left(\frac{\Ec(\theta) -
\Ec_{n}(\theta)}{|\theta-\theta^{*}|_{1}+\tilde{r}}\right)(|\Delta|_{1}+\tilde{r}).
\end{eqnarray}
By Lemma \ref{lem1} it holds on an event $\A_{1}$ of probability
at least $1-M^{-K}-3M^{-2K}\log\frac{n}{\log M}$ that
\begin{equation}\label{interm2}
\sup_{\theta \in \Theta : \theta\neq \theta^{*} }
\frac{\Ec(\theta) -
\Ec_{n}(\theta)}{|\theta-\theta^{*}|_{1}+\tilde{r}} \leqslant 2Kr.
\end{equation}

For any $\hat{\theta}\in\hat{\Theta}$, we have by
definition of the Dantzig estimator that $|\hat{\theta}|_{1}
\leqslant |\theta^{*}|_{1}$. Thus
\begin{equation}\label{Dsparse}
|\Delta_{J(\theta^{*})^{c}}|_{1} = \sum_{j\in
J(\theta^{*})^{c}}|\hat{\theta}_{j}| \leqslant \sum_{j\in
J(\theta^{*})}|\theta_{j}^{*}|-|\hat{\theta}_{j}|\leqslant|\Delta_{J(\theta^{*})}|_{1}.
\end{equation}

Define the function $g\,:\,t\rightarrow
R_{n}(\theta^{*}+t\Delta)$. Clearly $g$ is convex and
differentiable on $[0,1]$. Thus, the function $g'$ is
nondecreasing on $[0,1]$ with derivative $g'(t)=\nabla
R_{n}(\theta^{*} + t\Delta)^{T}\Delta$. The constraint
$\left|\nabla \hat{R}_{n}(\theta)\right|_{\infty}\leqslant r$ in
(\ref{Dantso}) and Lemma \ref{lem3} yield, on an event $\A_{2}$
of probability at least $1-M^{-1}$,
\begin{eqnarray}\label{interm3}
\Ec_{n}(\hat{\theta})&=&R_{n}(\hat{\theta})-R_{n}(\theta^{*})\nonumber\\
&=&\int_{0}^{1}\nabla
R_{n}(\theta^{*}+t\Delta)^{T}\Delta dt\nonumber\\
&\leqslant& r|\Delta|_{1},
\end{eqnarray}
for some numerical constant $C>0$.

Combining (\ref{interm1})-(\ref{interm3}) yields that on the event
$\A_{1}\cap \A_{2}$
\begin{equation}\label{interm4}
\Ec(\hat{\theta}) \leqslant (2+4K)r|\Delta_{J(\theta^{*})}|_{1}+
12\|\partial_{u}\gamma\|_{\infty}K\tilde{r}^{2}.
\end{equation}

Next,
\begin{eqnarray}\label{interm5}
2(1+2K)r|\Delta_{J(\theta^{*})}|_{1} &\leqslant& 2(1+2K)r\sqrt{s}|\Delta_{J(\theta^{*})}|_{2}\nonumber\\
&\leqslant&
\frac{2(1+2K)cr\sqrt{s}}{\zeta(s)}\frac{\|f_{\Delta}\|}{c}\nonumber\\
&\leqslant& \frac{1}{\kappa'}\left( \frac{2cr\sqrt{s}}{ \zeta(s)}
\right)^{\kappa'}+
\frac{1}{\kappa}\left(\frac{\|f_{\Delta}\|}{c}\right)^{\kappa}\nonumber\\
&\leqslant& \frac{1}{\kappa'}\left( \frac{2(1+2K)cr\sqrt{s}}{
\zeta(s)} \right)^{\kappa'}+ \frac{1}{\kappa}
\Ec(\hat{\theta}^{D}),
\end{eqnarray}
where we have used the Cauchy-Schwarz inequality in the first
line, the inequality $xy\leqslant
|x|^{\kappa}/\kappa+|y|^{\kappa'}/\kappa'$ that holds for any
$x,y$ in $\R$ and for any $\kappa,\kappa'$ in $(1,\infty)$ such
that $1/\kappa + 1/\kappa'=1$ in the third line, and Assumption 2
in the last line. Combining (\ref{interm4}) and (\ref{interm5})
and the fact that $\tilde{r}\leqslant 1$ yields the first
inequality. The second inequality is a consequence of
(\ref{th1-pred}) and (\ref{interm5}).
\end{proof}

We state and prove below intermediate results used in the proof of
Theorem 1.

\begin{lem}\label{lem1}
Let Assumptions \ref{A?-norm-sup} and \ref{A?-techass} be
satisfied. Then, with probability at least
$1-M^{-K}-3M^{-2K}\log\frac{n}{\log M}$, we have
\begin{equation}
\sup_{\theta \in \Theta} \frac{|\Ec(\theta) -
\Ec_{n}(\theta)|}{|\theta-\theta^{*}|_{1}+\tilde{r}} \leqslant
2Kr,
\end{equation}
where $r$ is defined in Theorem 1.
\end{lem}

\begin{proof}
For any $A>0$, define the random variable
\begin{equation*}
T_{A} = \sup_{\theta \in \Theta : |\theta-\theta^{*}|_{1}
\leqslant A}|\Ec_{n}(\theta)-\Ec(\theta)|.
\end{equation*}
For any $\theta$ in $\Theta$ and $(x,y)$ in $\Z$ we have
\begin{equation*}
|\gamma(y,f_{\theta}(x)) - \gamma(y,f_{\theta^{*}}(x))| \leqslant
\|\partial_{u}\gamma\|_{\infty}\left(
L|\theta-\theta^{*}|_{1}\wedge 2K \right),
\end{equation*}
and
\begin{equation*}
\E\left(|\gamma(Y,f_{\theta}(X)) -
\gamma(Y,f_{\theta^{*}}(X))|^{2}\right) \leqslant
\|\partial_{u}\gamma\|_{\infty}^{2}\left(
|\theta-\theta^{*}|_{1}^{2}\wedge 2K^{2} \right).
\end{equation*}

Assumption 3 and Bousquet's concentration inequality (cf. Theorem
4 in Section 6 below) with $x=(A\vee 2K)\log M$,
$c=2\|\partial_{u}\gamma\|_{\infty}(AL\wedge 2K)$ and $\sigma =
\sqrt{2}\|\partial_{u}\gamma\|_{\infty}(A\wedge \sqrt{2}K)$ yield
\begin{equation*}
\Prob\left(T_{A}\geqslant \E(T_{A}) +
2AK\|\partial_{u}\gamma\|_{\infty}\tilde{r}\right)\leqslant
M^{-(2K)\vee A}.
\end{equation*}
We study now the quantity $\E(T_{A})$. By standard symmetrization
and contraction arguments (cf. Theorems 5 and 6 in Section 6) we
obtain
\begin{eqnarray*}
\E(T_{A}) &\leqslant& 4 \|\partial_{u}\gamma\|_{\infty} \E\left(
\sup_{\theta \in \Theta\,:\, |\theta-\theta^{*}|_{1}\leqslant A }
\left|
\frac{1}{n}\sum_{i=1}^{n}\epsilon_{i}f_{\theta-\theta^{*}}(X_{i})
\right| \right).
\end{eqnarray*}
Then, observe that the mapping $u \rightarrow
\frac{1}{n}\sum_{i=1}^{n}\epsilon_{i}f_{u}(X_{i})$ is linear,
thus its supremum on a simplex is attained at one of its
vertices. This yields
\begin{equation*}
\E(T_{A}) \leqslant 4 \|\partial_{u}\gamma\|_{\infty} A \E\left(
\max_{1\leqslant j \leqslant M } \left|
\frac{1}{n}\sum_{i=1}^{n}\epsilon_{i}f_{j}(X_{i}) \right| \right).
\end{equation*}
Combining Assumption \ref{A?-techass} and Lemma 2 we obtain
\begin{equation*}
\E(T_{A})\leqslant4 \|\partial_{u}\gamma\|_{\infty} A\tilde{r}.
\end{equation*}
Thus
\begin{equation}\label{intpeeling1}
\Prob\left( T_{A}\geqslant 6AK\|\partial_{u}\gamma\|_{\infty}
\tilde{r} \right) \leqslant M^{-(2K)\vee A}.
\end{equation}

Define the following subsets of $\Theta$
\begin{eqnarray*}
\Theta(I) &=& \left\{ \theta \in \Theta \, : \,
|\theta-\theta^{*}|_{1} \leqslant \tilde{r} \right\},\\
\Theta(II) &=& \left\{ \theta \in \Theta \, : \,
\tilde{r}< |\theta-\theta^{*}|_{1} \leqslant 2K \right\},\\
\Theta(III) &=& \left\{ \theta \in \Theta \, : \,
|\theta-\theta^{*}|_{1} > 2K \right\}.
\end{eqnarray*}
For any $t>0$ define the probabilities
\begin{eqnarray*}
P_{I} &=& \Prob\left(\sup_{\theta \in \Theta(I)}
\frac{|\Ec(\theta) -
\Ec_{n}(\theta)|}{|\theta-\theta^{*}|_{1}+\tilde{r}}
\geqslant t\right)\\
P_{II} &=& \Prob\left(\sup_{\theta \in \Theta(II)}
\frac{|\Ec(\theta) -
\Ec_{n}(\theta)|}{|\theta-\theta^{*}|_{1}+\tilde{r}} \geqslant
t\right)\\
P_{III} &=& \Prob\left(\sup_{\theta \in \Theta(III)}
\frac{|\Ec(\theta) -
\Ec_{n}(\theta)|}{|\theta-\theta^{*}|_{1}+\tilde{r}} \geqslant
t\right)
\end{eqnarray*}
For any $t>0$ we have
\begin{align*}
\Prob\left(\sup_{\theta \in \Theta} \frac{|\Ec(\theta) -
\Ec_{n}(\theta)}{|\theta-\theta^{*}|_{1}+\tilde{r}} \geqslant
t\right) \leqslant P_{I}+P_{II}+P_{III}.
\end{align*}
Now, we bound from above the three probabilities on the right hand
side of the above expression. Take
$t=12\|\partial_{u}\gamma\|_{\infty}K\tilde{r}$. Applying
(\ref{intpeeling1}) to $P_{I}$ yields that
\begin{equation*}
P_{I}\leqslant \Prob\left( T_{\tilde{r}}\geqslant
6\|\partial_{u}\gamma\|_{\infty}K\tilde{r}^{2}\right)\leqslant
M^{-2K},
\end{equation*}
since we have $\tilde{r}\leqslant K$ by Assumption
\ref{A?-techass}.

Consider now $P_{II}$. We have
\begin{equation*}
\Theta(II)\subset  \bigcup_{j=0}^{j_{0}}\left\{ \theta\in\Theta\,
: \, A_{j+1}\leqslant |\theta-\theta^{*}|_{1}\leqslant A_{j}
\right\},
\end{equation*}
where $A_{j}=2^{1-j}K$, $j=0,\ldots,j_{0}$ and $j_{0}$ is chosen
such that $2^{1-j_{0}}K>\tilde{r} $ and $2^{-j_{0}}K\leqslant
\tilde{r} $. Thus
\begin{eqnarray*}
P_{II}&\leqslant& \sum_{j=0}^{j_{0}}\Prob\left(T_{A_{j}}\geqslant
12\|\partial_{u}\gamma\|_{\infty}A_{j+1}K\tilde{r}\right)\\
&\leqslant& \sum_{j=0}^{j_{0}}\Prob\left(T_{A_{j}}\geqslant
6\|\partial_{u}\gamma\|_{\infty}A_{j}K\tilde{r}\right)\\
&\leqslant& (j_{0}+1)M^{-2K}\\
&\leqslant& \left(3\left(\log\frac{n}{\log M}\right)-1
\right)M^{-2K}.
\end{eqnarray*}
Consider finally $P_{III}$. We have
\begin{equation*}
\Theta(III)\subset \bigcup_{j=0}^{\infty}\left\{
\theta\in\Theta\, : \, \bar{A}_{j-1} <
|\theta-\theta^{*}|_{1}\leqslant \bar{A}_{j}\right\},
\end{equation*}
where $\bar{A}_{j} = 2^{1+j}K$, $j\geqslant 0$. Thus
\begin{eqnarray*}
P_{III}&\leqslant&
\sum_{j=1}^{\infty}\Prob\left(T_{\bar{A}_{j}}\geqslant
12\|\partial_{u}\gamma\|_{\infty}\bar{A}_{j-1}K\tilde{r}\right)\\
&\leqslant& \sum_{j=0}^{j_{0}}\Prob\left(T_{A_{j}}\geqslant
6\|\partial_{u}\gamma\|_{\infty}\bar{A}_{j}K\tilde{r}\right)\\
&\leqslant& \sum_{j=1}^{\infty}M^{-\bar{A}_{j}}\\
&\leqslant& M^{-K}.
\end{eqnarray*}

\end{proof}

We now study the quantity $ \E\left(\max_{1\leqslant j \leqslant
M} \left| \frac{1}{n}\sum_{i=1}^{n}\epsilon_{i}f_{j}(X_{i})
\right|\right)$. This is done in the next lemma.

\begin{lem}\label{lem2}
We have
\begin{equation}\label{Bernsteintype}
\E\left(\max_{1\leqslant j \leqslant M} \left|
\frac{1}{n}\sum_{i=1}^{n}\epsilon_{i}f_{j}(X_{i})
\right|\right)\leqslant  \tilde{r},
\end{equation}
where $\tilde{r}$ is defined in (\ref{rtilde}).
\end{lem}

\begin{proof}
Define the random variables
\begin{equation*}
U_{j}=\frac{1}{\sqrt{n}}\sum_{i=1}^{n}\epsilon_{i}f_{j}(X_{i}).
\end{equation*}
The Bernstein inequality yields, for any $j=1,\ldots,M$ and $t>0$,
\begin{equation}
\Prob\left( |U_{j}| \geqslant t \right) \leqslant \exp\left( -
\frac{t^{2}}{2(t\|f_{j}\|_{\infty}/(3\sqrt{n})+\|f_{j}\|^{2})}
\right).
\end{equation}

Set $b_{j}=\| f_{j} \|_{\infty}/(3\sqrt{n})$. Define the random
variables $T_{j}=U_{j}\1_{|Y_{j}|>\|f_{j}\|^{2}/b_{j}}$ and
$T'_{j}=U_{j}\1_{|Y_{j}|\leqslant \|f_{j}\|^{2}/b_{j}}$. For all
$ t>0$ we have
\begin{equation*}
\Prob\left( |T_{j}|>t\right)\leqslant 2\exp\left(-
\frac{t}{4b_{j}} \right),\hspace{0.25cm} \Prob\left(
|T'_{j}|>t\right)\leqslant 2\exp\left(-
\frac{t^{2}}{4\|f_{j}\|^{2}} \right).
\end{equation*}
Define the function $h_{\nu}(x)=\exp(x^{\nu})-1$, where $\nu>0$.
This function is clearly convex for any $\nu>0$. We have
\begin{equation*}
\E\left( h_{1}\left( \frac{|T_{j}|}{12b_{j}} \right) \right) =
\int_{0}^{\infty}e^{t}\Prob(|T_{j}|>12b_{j}t)dt\leqslant 1,
\end{equation*}
where we have used Fubini's Theorem in the first equality. Since
the function $h_{1}$ is convex and nonnegative, we have
\begin{eqnarray*}
h_{1}\left( \E\left(\max_{1\leqslant j \leqslant M
}\frac{|T_{j}|}{12b_{j}} \right) \right)&\leqslant&
\E\left(h_{1}\left(\max_{1\leqslant j \leqslant M
}\frac{|T_{j}|}{12b_{j}}\right) \right)\\
&\leqslant&
\E\left(\sum_{j=1}^{M}h_{1}\left(\frac{|T_{j}|}{12b_{j}}
\right)\right)\\
&\leqslant& M,
\end{eqnarray*}
where we have used the Jensen inequality. Since the function
$h_{1}^{-1}(x)=\log(1+x)$ is increasing, we have
\begin{eqnarray}\label{inter3}
\E\left(\max_{1\leqslant j \leqslant M }\frac{|T_{j}|}{12b_{j}}
\right)&\leqslant& \log(1+M)\nonumber\\
\E\left(\max_{1\leqslant j \leqslant M }|T_{j}|\right)&\leqslant&
4\frac{\log(1+M)}{\sqrt{n}}\max_{1\leqslant j \leqslant M
}\|f_{j}\|_{\infty}.
\end{eqnarray}
Applying the same argument to the function $h_{2}$, we prove that
\begin{equation}\label{inter4}
\E\left(\max_{1\leqslant j \leqslant M }|T_{j}'| \right)\leqslant
2\sqrt{3}\sqrt{\log(1+M)}\max_{1\leqslant j \leqslant M
}\|f_{j}\|.
\end{equation}
Combining (\ref{inter3}) and (\ref{inter4}) yields the result.

\end{proof}

\begin{lem}\label{lem3}
Let Assumptions \ref{A?-norm-sup} and \ref{A?-techass} be
satisfied. Then, with probability at least $1-M^{-1}$, we have
$$|\nabla
\hat{R}_{n}(\theta^{*})|_{\infty}\leqslant r,$$
where $r$ is
defined in Theorem 1.
\end{lem}

\begin{proof}
For any $1\leqslant j \leqslant M$ define
\begin{equation*}
Z_{j} =
\frac{1}{n}\sum_{i=1}^{n}\partial_{u}\gamma(Y_{i},f_{\theta^{*}}(X_{i}))f_{j}(X_{i}).
\end{equation*}
Since the function $\theta\rightarrow \gamma(y,f_{\theta}(x))$ is
differentiable w.r.t. $\theta$ and
$|\partial_{u}\gamma(y,f_{\theta}(x))f_{j}(x)|\leqslant
\|\partial_{u} \gamma\|_{\infty}L$ for any $(x,y)\in\X\times\Y$
and $\theta \in \R^{M}$, we have
\begin{equation*}
\E(Z_{j})=\frac{\partial R(\theta^{*})}{\partial\theta_{j}}=0.
\end{equation*}

Next, similarly as in Lemmas \ref{lem1} and \ref{lem2}, we prove
that
\begin{equation*}
\E(|\nabla \hat{R}_{n}(\theta^{*})|_{\infty})  \leqslant
4\|\partial_{u}\gamma\|_{\infty}\tilde{r}.
\end{equation*}
Finally Bousquet's concentration inequality (cf. Theorem 4 in
Section 6 below) yields that, with probability at least
$1-M^{-1}$,
\begin{align*}
|\nabla \hat{R}_{n}(\theta^{*})|_{\infty} &\leqslant
\E(|\nabla \hat{R}_{n}(\theta^{*})|_{\infty})\\
&\hspace{0.25cm}+
\sqrt{2\frac{\log M}{n}\left( \|\partial_{u}\gamma\|_{\infty}^{2}
+ 2\|\partial_{u}\gamma\|_{\infty}L
\E(|\nabla \hat{R}_{n}(\theta^{*})|_{\infty})\right)}\\
&\hspace{0.25cm}+\frac{\|\partial_{u}\gamma\|_{\infty}L \log M}{3n}\\
&\leqslant 6\|\partial_{u}\gamma\|_{\infty}\tilde{r}.
\end{align*}

\end{proof}

\section{Examples}

\subsection{Robust regression with the Huber loss}

We consider the linear regression model
\begin{equation}\label{lin-reg}
Y=f_{\theta^{*}}(X)+W,
\end{equation}
where $X\in \R^{d}$ is a random vector, $W\in\R$ is a random
variable independent of $X$ whose distribution is symmetric w.r.t.
$0$ and $\theta^{*}\in\R^{M}$ is the unknown vector of
parameters. Consider the function
\begin{equation*}
\phi(x) = \frac{1}{2}x^{2}\1_{|x|\leqslant 2K+\alpha} + \left(
(2K+\alpha)|x|-\frac{(2K+\alpha)^{2}}{2}
\right)\1_{|x|>2K+\alpha},
\end{equation*}
where $\alpha>0$. The Huber loss function is defined by
\begin{equation}\label{Huberloss}
Q(z,\theta) = \phi(y-f_{\theta}(x)),
\end{equation}
where $z=(x,y) \in \R^{d}\times \R$ and $\theta \in \Theta$.

In the following lemma we prove that for this loss function Assumption \ref{A2-MA} is
satisfied with $\kappa = 2$ and $c=(2/\Prob(|W| \leqslant
\alpha))^{1/2}$.
\begin{lem}\label{Cond-MA}
Let $Q$ be defined by (\ref{Huberloss}). Then for any $\theta\in\Theta$ we have
\begin{equation*}
\frac{\Prob(|W|\leqslant
\alpha)}{2}\|f_{\theta}-f_{\theta^{*}}\|^{2}\leqslant \Ec(\theta).
\end{equation*}
\end{lem}
\begin{proof}
Set $\Delta=\theta-\theta^{*}$. Since $\phi'$ is absolutely continuous, we have for any
$\theta\in\Theta$
\begin{align*}
Q(Z,\theta)-Q(Z,\theta^{*}) &= \phi'(W)f_{-\Delta}(X)\\
&\hspace{0.25cm}+
\left[\int_{0}^{1}\1_{|W+tf_{-\Delta}(X)|\leqslant 2K+\alpha}(1-t)dt\right]f_{\Delta}(X)^{2}\\
&\geqslant \phi'(W)f_{-\Delta}(X) + \frac{1}{2}\1_{(|W|\leqslant
\alpha)}f_{\Delta}(X)^{2},
\end{align*}
since
$\|f_{\theta}\|_{\infty}\leqslant K$ for any $\theta\in\Theta$.
Taking the expectations we get
\begin{equation*}
R(\theta)-R(\theta^{*}) \geqslant \frac{\Prob(|W|\leqslant
\alpha)}{2}\|f_{\Delta}\|^{2},
\end{equation*}
for any $\alpha>0$ since $\phi'$ is odd and the distribution of $W$
is symmetric w.r.t. $0$.
\end{proof}
We have the following corollary of Theorem 1.
\begin{cor}\label{cor1}
Let Assumptions \ref{A1-RE}, \ref{A?-norm-sup} and
\ref{A?-techass} be satisfied. If $M(\theta^{*})\leqslant s$,
then, with probability at least $1-M^{-1}-M^{-K}-3M^{-2K}\log
\frac{n}{\log M}$, we have
\begin{equation*}
\sup_{\hat{\theta}\in\hat{\Theta}}\Ec(\hat{\theta})\leqslant
\frac{8(1+2K)^{2}}{\Prob(|W|\leqslant \alpha)\zeta(s)^{2}}sr^{2}
+ \frac{2}{3}r^{2},
\end{equation*}
and
\begin{equation*}
\sup_{\hat{\theta}\in\hat{\Theta}}|\hat{\theta}-\theta^{*}|_{1}\leqslant
\frac{8(1+2K)}{\Prob(|W|\leqslant \alpha)\zeta(s)^{2}}sr
+\frac{K}{3(1+2K)}r.
\end{equation*}
\end{cor}

\subsubsection{Logistic regression and similar models}
We consider $Z=(X,Y)\in \X\times \{0,1\}$ where $\X$ is a Borel subset of
$\R^{d}$. The conditional probability $\Prob(Y=1\,|X=x)=\pi(x)$
is unknown where $\pi$ is a function on $\X$ with values in
$[0,1]$. We assume that $\pi$ is of the form
\begin{equation}\label{logit}
\pi(x)=\Phi'(f_{\theta^{*}}(x)),
\end{equation}
where the function $\Phi\, : \, \R \rightarrow \R^{*}$ is convex,
twice differentiable, of derivative $\Phi'$ with values in $[0,1]$
and the vector $\theta^{*}\in\R^{M}$ is unknown. Consider, e.g.,
the logit loss function $\Phi(u)=\log(1+e^{u})$. We assume that
$\Phi$ is known. Define the quantity
\begin{equation}\label{tau} \tau(R)=\frac{1}{2}\inf_{|u|\leqslant
R}\Phi^{(2)}(u),
\end{equation}
for any $R\geqslant 0$. We want to estimate $\theta^{*}$ with the
procedure (\ref{Dantso}) and the convex loss function
\begin{equation}\label{logitmodel}
Q(z,\theta)=-yf_{\theta}(x)+\Phi(f_{\theta}(x)),
\end{equation}
where
$z=(x,y)\in \R^{d}\times \{0,1\}$. Thus we need to check
Assumption \ref{A2-MA} to apply Theorem 1.
\begin{lem}
Let the loss function be of the form (\ref{logitmodel}) where $\Phi$ satisfies the above assumptions. Then for any $\theta\in\R^{M}$ we have
\begin{equation*}
\tau(K)\| f_{\theta}-f_{\theta^{*}}\|^{2} \leqslant \Ec(\theta).
\end{equation*}
\end{lem}

\begin{proof}
For any $\theta \in\Theta$, we have
\begin{align*}
Q(Z,\theta)-Q(Z,\theta^{*})&=
\nabla Q(Z,\theta^{*})^{T}(\theta-\theta^{*})\\
&\hspace{0.25cm}+
\left[\int_{0}^{1}\Phi^{(2)}(H(X)^{T}(\theta^{*}+t(\theta-\theta^{*})))(1-t) dt\right] f_{\Delta}(X)^{2}\\
&\geqslant \nabla Q(Z,\theta^{*})^{T}(\theta-\theta^{*}) +
\tau(\|f_{\theta}\|_{\infty}\vee \|f_{\theta^{*}}\|_{\infty})
 f_{\Delta}(X)^{2}.
\end{align*}
Since $\|\nabla Q(\cdot,\cdot)\|_{\infty}<\infty$, we can
differentiate under the expectation sign, so that
\begin{equation*}
\E(\nabla Q(Z,\theta^{*})^{T}(\theta-\theta^{*}))=\nabla
R(\theta^{*})=0.
\end{equation*}
Thus
\begin{equation*}
\Ec(\theta) \geqslant \tau(\|f_{\theta}\|_{\infty}\vee
\|f_{\theta^{*}}\|_{\infty})\|f_{\theta}-f_{\theta^{*}}\|^{2}.
\end{equation*}
\end{proof}
Thus Assumption \ref{A2-MA} is satisfied with the constants
$\kappa=2$ and $c=\frac{1}{ \sqrt{\tau(K)}}$. We have the
following corollary of Theorem 1.
\begin{cor}
Let Assumptions \ref{A1-RE}, \ref{A?-norm-sup} and
\ref{A?-techass} be satisfied. If $M(\theta^{*})\leqslant s$,
then, with probability at least $1-M^{-1}-M^{-K}-3M^{-2K}\log
\frac{n}{\log M}$, we have
\begin{equation*}
\sup_{\hat{\theta}\in\hat{\Theta}}\Ec(\hat{\theta})\leqslant
\frac{4(1+2K)^{2}}{\tau(K)\zeta(s)^{2}}sr^{2} + \frac{2}{3}r^{2},
\end{equation*}
and
\begin{equation*}
\sup_{\hat{\theta}\in\hat{\Theta}}|\hat{\theta}-\theta^{*}|_{1}\leqslant
\frac{4(1+2K)}{\tau(K)\zeta(s)^{2}}sr + \frac{K}{3(1+2K)}r.
\end{equation*}
\end{cor}

\section{Sup-norm convergence rate for the regression model with the Huber loss}
In this section, we derive the sup-norm convergence rate of the
Dantzig estimators to the target vector $\theta^{*}$ in the
linear regression model under a mutual coherence assumption on the
dictionary and Huber's loss. The proof relies on the fact that the Hessian matrix
of the risk also satisfies the mutual coherence condition for
this particular model. Unfortunately, we cannot proceed similarly
in the general case because the Hessian matrix of the risk at
point $\theta^{*}$ does not necessarily satisfy the mutual
coherence condition even if the Gram matrix of the dictionary
does. Note that for Huber's loss the Dantzig minimization problem (\ref{Dantso}) is computable feasible. The constraints in (\ref{Dantso}) are indeed linear, so that (\ref{Dantso}) is a linear programming problem.

Denote by $\Psi(\theta)$ the Hessian matrix of the risk $R$
evaluated at $\theta$. With our assumptions on the dictionary
$\D$ and on the function $\gamma$, for any $\theta\in\R^{M}$
we have
\begin{equation*}
\Psi(\theta) \stackrel{\triangle}{=} \nabla^{2} R(\theta) =
\left(\E\left(
\partial^{2}_{u}\gamma(Y,f_{\theta}(X))f_{j}(X)f_{k}(X) \right)\right)_{1\leqslant
j,k \leqslant M}.
\end{equation*}
Note that for the quadratic loss we have $\Psi(\cdot)\equiv 2G$
where $G$ is the Gram matrix of the design. For Lipschtiz loss
functions the Hessian matrix $\Psi$ varies with $\theta$.

We consider the linear regression model (\ref{lin-reg}). For any
functions $g,h\, : \X\rightarrow \R$, denote by $<g,h>$ the
scalar product $\E(g(X)h(X))$. Define the Gram matrix $G$ by
\begin{equation*}
G=(<f_{j},f_{k}>)_{1\leqslant j, k \leqslant M}.
\end{equation*}
From now on, we assume that $G$ satisfies a mutual coherence
condition.
\begin{assum}\label{A3-coh}
The Gram matrix $G=(<f_{j},f_{k}>)_{1\leqslant j,k\leqslant M }$
satisfies
\begin{equation*}
G_{j,j}=1, \, \forall 1\leqslant j \leqslant M,
\end{equation*}
and
\begin{equation*}
\max_{j \neq k}|G_{j,k}| \leqslant \frac{1}{3 \beta s},
\end{equation*}
where $s\geqslant 1$ is an integer and $\beta > 1$ is a constant.
\end{assum}
This assumption is stronger than Assumption \ref{A1-RE}. We have
indeed the following Lemma (cf. Lemma 2 in \cite{L08}).
\begin{lem}\label{coh-RE}
Let Assumption \ref{A3-coh} be satisfied. Then
\begin{equation*}
\zeta(s)\stackrel{\triangle}{=}\min_{J\subset \{1,\cdots,M
\},|J|\leqslant s}\,\min_{\Delta\neq 0:|\Delta_{J^{c}}|_{1}\leqslant
|\Delta_{J}|_{1}}\frac{\|f_{\Delta}\|}{|\Delta_{J}|_{2}}\geqslant
\sqrt{1-\frac{1}{\beta}}>0.
\end{equation*}
\end{lem}

Note that Assumption \ref{A3-coh} the
vector $\theta^{*}$ satisfying (\ref{lin-reg}) such that
$M(\theta^{*})\leqslant s$ is \textbf{unique}. Consider indeed two
vectors $\theta^{1}$ and $\theta^{2}$ satisfying (\ref{lin-reg})
such that $M(\theta^{1})\leqslant s$ and $M(\theta^{2})\leqslant
s$. Denote $\theta=\theta^{1}-\theta^{2}$ and $J=J(\theta^{1})\cup
J(\theta^{2})$. Clearly we have $f_{\theta}(X)=0$ a.s. and
$M(\theta)\leqslant 2s$. Assume that $\theta^{1}$ and $\theta^{2}$
are distinct. Then,
\begin{eqnarray*}
\frac{\|f_{\theta}\|_{2}^{2}}{|\theta|_{2}^{2}} &=& 1 +
\frac{\theta^{T}(G-I_{M})\theta}{|\theta|_{2}^{2}}\\
&\geqslant& 1-\frac{1}{3\beta
s}\sum_{i,j=1}^{M}\frac{|\theta_{i}||\theta_{j}|}{|\theta|_{2}^{2}}\\
&\geqslant& 1-\frac{1}{3\beta}>0,
\end{eqnarray*}
where we have used the Cauchy-Schwarz inequality. This contradicts
the fact that $f_{\theta}(X)=0$ a.s.

For the linear regression model, the Hessian matrix $\Psi$ at point
$\theta$ is
\begin{equation*}
\Psi(\theta) = \E(\1_{|f_{\theta^{*}-\theta}(X)+W|\leqslant
2K+\alpha }f_{j}(X)f_{k}(X))_{1\leqslant j,k \leqslant M}.
\end{equation*}
We observe that
\begin{equation*}
\Psi(\theta^{*})=\Prob(|W|\leqslant 2K+\alpha)G.
\end{equation*}
Thus $\Psi(\theta^{*})$ satisfies a condition similar to Assumption 4 but with a different constant if
$\Prob(|W|\leqslant 2K+\alpha)>0$. The empirical Hessian matrix
$\hat{\Psi}$ at point $\theta\in\R^{M}$ is defined by
\begin{equation*}
\hat{\Psi}_{j,k}(\theta) =
\frac{1}{n}\sum_{i=1}^{n}\1_{|f_{\theta^{*}-\theta}(X_{i})+W_{i}|\leqslant
2K+\alpha }f_{j}(X_{i})f_{k}(X_{i}),\hspace{0.1cm}1\leqslant j,k
\leqslant M.
\end{equation*}

We will prove that the empirical Hessian matrix
$\hat{\Psi}(\theta)$ satisfies a mutual coherence condition for
any $\theta$ in a small neighborhood of $\theta^{*}$ under some
additional assumptions given below.

First, we need an additional mild assumption on the noise.
\begin{assum}\label{A6-noise}
There c.d.f. $F_{W}$ of $W$ is Lipschitz continuous.
\end{assum}
This assumption is satisfied, e.g., if $W$ admits a bounded
density so we allow heavy tailed distributions such as the Cauchy. In the
sequel we assume w.l.o.g. that the Lipschitz constant of $F_{W}$
equals $1$.

We impose a restriction on the sparsity $s$.
\begin{assum}\label{A5-sparsity}
The sparsity $s$ satisfies $s\leqslant \frac{1}{\sqrt{r}}$.
\end{assum}
This implies that we can recover the sparse vectors with at most
$O\left(\left(n/\log M\right)^{1/4}\right)$ nonzero components.

Define $V_{\eta} = \{ \theta \in \Theta\, : \,
|\theta-\theta^{*}|_{1}\leqslant \eta \}$ where $\eta = C_{1}rs$
and
\begin{equation}\label{const-C1}
C_{1} = \frac{8(1+2K)\beta}{\Prob(|W|\leqslant \alpha)(\beta-1)}
+ \frac{1}{6}.
\end{equation}
Consider the event
\begin{equation}\label{Evt1}
E=\left\{ \sup_{ 1\leqslant j, k \leqslant M, \theta \in V_{\eta}
}\left| \hat{\Psi}_{j,k}(\theta)- \Psi_{j,k}(\theta) \right|
\leqslant 8L^{3}\eta + 4L\tilde{r} +  \frac{C_{2}}{\sqrt{n}s^{2}}
\right\},
\end{equation}
where
$$C_{2} =  2\sqrt{1+(1+L^{2})\left( 8C_{1}L^{3}+ \frac{4L}{s}\right)} + \frac{1+L^{2}}{3}.$$

We have the following intermediate result.
\begin{lem}\label{lem7}
Let Assumptions \ref{A?-norm-sup}- \ref{A6-noise} be satisfied.
Then $\Prob(E) \geqslant 1-\exp(-\sqrt{\log M})$.
\end{lem}

\begin{proof}
Define the variable
\begin{equation*}
Z=\sup_{1\leqslant j, k\leqslant M,\, \theta \in V_{\eta}
}\left|\hat{\Psi}_{j,k}(\theta)-\Psi_{j,k}(\theta)\right|.
\end{equation*}
Applying the Bousquet concentration inequality (cf. Theorem 4 in Section 6) with the constants
$c=(1+L^{2})/n$, $\sigma^{2}=2/n^{2}$ and $x=\frac{\sqrt{n}}{s^{2}}$ yields that, with probability at least
$1-e^{-x}$,
\begin{eqnarray}\label{interm6}
Z&\leqslant&
\E\left(Z\right)+\frac{2}{\sqrt{n}s}\sqrt{1+(1+L^{2})\E\left(Z\right)}+\frac{1+L^{2}}{3\sqrt{n}s^{2}}.
\end{eqnarray}
We study now the quantity $\E(Z)$. A standard symmetrization and
contraction argument yields
\begin{align}\label{interm7}
\E(Z) &\leqslant 2\E\left( \sup_{1\leqslant j,k\leqslant M,\,
\theta\in
V_{\eta}}\left|\frac{1}{n}\sum_{i=1}^{n}\epsilon_{i}\1_{|f_{\theta^{*}-\theta}(X_{i})+W_{i}|\leqslant
2K+\alpha}f_{j}(X_{i})f_{k}(X_{i})\right| \right)\nonumber\\
&\leqslant
2\E\left(\left|\frac{1}{n}\sum_{i=1}^{n}\epsilon_{i}\1_{|W_{i}|\leqslant
2K+\alpha}f_{j}(X_{i})f_{k}(X_{i})\right| \right)\nonumber\\
&+ 2\E\left(\sup_{1\leqslant j,k\leqslant M,\, \theta\in
V_{\eta}}\left|\frac{1}{n}\sum_{i=1}^{n}\epsilon_{i}(\1_{|f_{\theta^{*}-\theta}(X_{i})+W_{i}|\leqslant
2K+\alpha}-\1_{|W_{i}|\leqslant
2K+\alpha})f_{j}(X_{i})f_{k}(X_{i})\right| \right).
\end{align}
Denote by $(I)$ and $(II)$ respectively the first term and the
second term on the right hand side of the above expression. The
contraction principle yields
\begin{equation}\label{interm8}
(I) \leqslant 4 \E\left( \max_{1\leqslant j,k \leqslant M}\left|
\frac{1}{n}\sum_{i=1}^{n}\epsilon_{i}f_{j}(X_{i})f_{k}(X_{i})
\right| \right).
\end{equation}
Then, similarly as in the proof of Lemma 2 we get
\begin{equation*}
\E\left( \max_{1\leqslant j,k \leqslant M}\left|
\frac{1}{n}\sum_{i=1}^{n}\epsilon_{i}f_{j}(X_{i})f_{k}(X_{i})
\right| \right) \leqslant L\tilde{r}.
\end{equation*}
Thus, for (II) we have
\begin{eqnarray}\label{interm9}
(II) &\leqslant& 2L^{2}\E\left( \sup_{ \theta \in V_{\eta}}
\frac{1}{n}\sum_{i=1}^{n}|\1_{|f_{\Delta}(X_{i})+W_{i}|\leqslant
2K+\alpha}-\1_{|W_{i}|\leqslant 2K+\alpha}|   \right)\nonumber\\
&\leqslant& 2L^{2}\Prob\left( 2K+\alpha-L\eta \leqslant
|W|\leqslant
2K+\alpha+L\eta\right)\nonumber\\
&\leqslant& 8L^{3}\eta.
\end{eqnarray}
Assumptions \ref{A?-techass} and \ref{A5-sparsity} yield that
$s\leqslant \left(\frac{n}{\log M}\right)^{1/4}$. Combining
(\ref{interm6})-(\ref{interm9}) yields the result.
\end{proof}

We need an additional technical assumption.
\begin{assum}\label{A?-techass2}
We have
$12L^{3}\eta+L\tilde{r}+\frac{C_{2}}{\sqrt{n}s^{2}}\leqslant
\frac{\Prob(|W|\leqslant 2K+\alpha)}{2}$.
\end{assum}
This is a mild assumption. It is indeed satisfied for $n$ large
enough if we assume that $\Prob(|W|\leqslant 2K+\alpha)>0$ since
Assumption 6 implies that $r\rightarrow 0$ as
$n\rightarrow \infty$.

We have the following result on the empirical Hessian matrix.
\begin{lem}\label{lem8}
Let Assumptions \ref{A?-norm-sup}-\ref{A?-techass2} be satisfied.
Then, with probability at least $1-\exp(-\sqrt{\log M})$, for any
$\theta\in V_{\eta}$, we have
\begin{eqnarray}\label{emp-coh}
\min_{1\leqslant j \leqslant
M}|\hat{\Psi}_{j,j}(\theta)|&\geqslant& \frac{\Prob(|W|\leqslant 2K+\alpha)}{2},\nonumber\\
\max_{ j\neq k } |\hat{\Psi}_{j,k}(\theta)|&\leqslant&
\frac{C_{3}}{s},
\end{eqnarray}
where $C_{3} = \frac{1}{3\beta} + 12 L^{3}C_{1} +
\frac{C_{2}}{\sqrt{n}s}$.
\end{lem}

\begin{proof}
For any $\theta$ in $V_{\eta}$ and any $j,k$ in $\{1,\ldots,M\}$ we
have
\begin{equation*}
\Psi_{j,k}(\theta)-\Psi_{j,k}(\theta^{*}) = \E\left(
(\1_{|f_{\Delta}(X)+W|\leqslant 2K+\alpha}-\1_{|W|\leqslant
2K+\alpha})f_{j}(X)f_{k}(X)   \right),
\end{equation*}
where $\Delta=\theta-\theta^{*}$. Then
\begin{eqnarray*}
|\Psi_{j,k}(\theta)-\Psi_{j,k}(\theta^{*})| &\leqslant&
L^{2}\E\left( |\1_{|f_{\Delta}(X)+W|\leqslant
2K+\alpha}-\1_{|W|\leqslant
2K+\alpha}|  \right)\\
&\leqslant& L^{2} \Prob\left(|W|\leqslant 2K+\alpha\, ,\,
|f_{\Delta}(X)+W|> 2K+\alpha \right)\\
&+& L^{2}\Prob\left(|W|> 2K+\alpha\, ,\, |f_{\Delta}(X)+W|\leqslant
2K+\alpha \right).
\end{eqnarray*}
Recall that $|f_{\Delta}(X)|\leqslant L\eta$. Then
\begin{eqnarray}\label{interm10}
|\Psi_{j,k}(\theta)-\Psi_{j,k}(\theta^{*})|&\leqslant&
L^{2}\Prob\left( 2K+\alpha -L\eta \leqslant |W| \leqslant 2K+\alpha
+ L\eta \right)\nonumber\\
&\leqslant& 2L^{2}\Prob\left( 2K+\alpha -L\eta \leqslant W \leqslant
2K+\alpha + L\eta \right)\nonumber\\
&\leqslant& 4L^{3}\eta,
\end{eqnarray}
where we have used the fact that the distribution of $W$ is
symmetric w.r.t. $0$ in the second line and Assumption
\ref{A6-noise} in the last line. Lemma \ref{lem7} and
(\ref{interm10}) yield that, on the event $E$, for any $\theta\in
V_{\eta}$,
\begin{equation*}
\min_{1\leqslant j \leqslant M} \hat{\Psi}_{j,j}(\theta) \geqslant
\Prob(|W|\leqslant
2K+\alpha)-12L^{3}\eta-\frac{C_{2}}{\sqrt{n}s^{2}},
\end{equation*}
and
\begin{equation*}
\max_{j\neq k}|\Psi_{j,k}(\theta)|\leqslant \frac{C_{3}}{s}.
\end{equation*}
\end{proof}

Now we can derive the optimal sup-norm convergence rate of the Dantzig estimators.
\begin{theo}\label{th2-supnorm}
Let Assumptions \ref{A?-norm-sup}-\ref{A?-techass2} be satisfied.
If $M(\theta^{*})\leqslant s$, then, on an event of probability at
least $1- M^{-1}-M^{-K} - \exp(-\sqrt{\log M}) -3M^{-2K}\log
\frac{n}{\log M}$, we have
\begin{equation*}
\sup_{\hat{\theta}\in
\hat{\Theta}}|\hat{\theta}-\theta^{*}|_{\infty} \leqslant C_{4}r,
\end{equation*}
where $r$ is defined in
Theorem 1,
\begin{equation*}
C_{4}=\frac{4+2C_{1}C_{3}}{\Prob(|W|\leqslant 2K+\alpha)},
\end{equation*}
with $C_{1}$ and $C_{3}$ defined respectively in (\ref{const-C1})
and Lemma \ref{lem8}.
\end{theo}

\begin{proof}
For any $\hat{\theta}$ in $\hat{\Theta}$ we have
\begin{equation*}
\nabla R_{n}(\hat{\theta})-\nabla R_{n}(\theta^{*}) =
\left[\int_{0}^{1}\hat{\Psi}(\theta^{*}+t\Delta)dt\right]\Delta,
\end{equation*}
where $\Delta=\hat{\theta}-\theta^{*}$.

The definition of our estimator, Lemma \ref{lem3} and
Corollary \ref{cor1} yield that, on an event $\A_{1}$ of
probability at least $1- M^{-1} - \exp(-\sqrt{\log M}) -3M^{-2K}\log \frac{n}{\log M}
$, we have that $\hat{\theta}\in
V_{\eta}$ and
\begin{equation*}
\left|
\left[\int_{0}^{1}\hat{\Psi}(\theta^{*}+t\Delta)dt\right]\Delta
\right|_{\infty} \leqslant 2r.
\end{equation*}
Lemma \ref{lem8} yields that, on the event $\A_{1} \cap E$,
\begin{equation*}
\frac{\Prob(|W|\leqslant 2K+\alpha)}{2}\left|\Delta
\right|_{\infty} \leqslant
2r + \frac{C_{3}}{s}|\Delta|_{1},
\end{equation*}
so that
\begin{equation*}
\left|\Delta \right|_{\infty} \leqslant C_{4}r.
\end{equation*}
\end{proof}
Note that Theorem 2 holds true for the Lasso estimators (2) with
exactly the same proof, provided that a result similar to Theorem 1 is valid for the Lasso estimators. This is in fact the case (cf.
\cite{VdG07,Kolt07}).

\section{Sign concentration property}

Now we study the sign concentration property of the Dantzig
estimators. We need an additional assumption on the magnitude of
the nonzero components of $\theta^{*}$.
\begin{assum}\label{A7-seuil}
We have
\begin{equation*}
\rho\stackrel{\Delta}{=}\min_{j\in J(\theta^{*})}|\theta^{*}_{j}|
> 2C_{4}r,
\end{equation*}
where $r$ is defined in Theorem 1 and $C_{4}$ is defined in
Theorem \ref{th2-supnorm}.
\end{assum}
We can find similar assumptions on $\rho$ in the work on sign
consistency of the Lasso estimator mentioned above. More
precisely, the lower bound on $\rho$ is of the order
$(s(\log M)/n)^{1/4}$ in \cite{MY06}, $n^{-\delta/2}$ with
$0<\delta<1$ in \cite{W06,ZY06}, $\sqrt{(\log Mn)/n}$ in
\cite{B07}, $\sqrt{s(\log M)/n}$ in \cite{ZH07} and $r$ in \cite{L08}.

We introduce the following thresholded version of our estimator. For any
$\hat{\theta}\in\hat{\Theta}$ the associated thresholded
estimator $\tilde{\theta}\in\R^{M}$ is defined by
\begin{equation*}
\tilde{\theta}_{j}=
\begin{cases}
\hat{\theta}_{j}, &\text{if
$|\hat{\theta}_{j}|>C_{4}r$},\\
0 &\text{elsewhere}.
\end{cases}
\end{equation*}
Denote by $\tilde{\Theta}$ the set of all such $\tilde{\theta}$.
We have first the following non-asymptotic result that we call
sign concentration property.

\begin{theo}\label{Signnonasymp}
Let Assumptions \ref{A?-norm-sup} and \ref{A3-coh}-\ref{A7-seuil}
be satisfied. If $M(\theta^{*})\leqslant s$, then
\begin{align*}
\Prob\left(\vec{\mathrm{sign}}(\tilde{\theta})=\vec{\mathrm{sign}}(\theta^{*}),\hspace{0.1cm}\forall
\tilde{\theta}\in\tilde{\Theta}\right)&\geqslant
1- M^{-1}-M^{-K}-
\exp(-\sqrt{\log M})\\
&\hspace{0.25cm}  -3M^{-2K}\log \frac{n}{\log M}.
\end{align*}
\end{theo}

Theorem \ref{Signnonasymp} guarantees that the sign vector of every vector
$\tilde{\theta}\in\tilde{\Theta}$ coincides with that of $\theta^{*}$ with probability close to one.

\begin{proof}
Theorem \ref{th2-supnorm} yields
$\sup_{\hat{\theta}\in\hat{\Theta}}|\hat{\theta}-\theta^{*}|_{\infty}\leqslant
C_{3}r$ on an event $\A$ of probability at least $1-6M^{-1}$. Take
$\hat{\theta}\in\hat{\Theta}$. For $j\in J(\theta^{*})^{c}$, we
have $\theta^{*}_{j}=0$, and $|\hat{\theta}_{j}|\leqslant c_{2}r$
on $\A$. For $j\in J(\theta^{*})$, we have
$|\theta^{*}_{j}|\geqslant 2C_{3}r$ and
$|\theta^{*}_{j}|-|\hat{\theta}_{j}|\leqslant
|\theta^{*}_{j}-\hat{\theta}C_{j}| \leqslant  C_{3}r$ on $\A$.
Since we assume that $\rho>2C_{3}$, we have on $\A$ that
$|\hat{\theta}_{j}|\geqslant >c_{2}r$. Thus on the event $\A$ we
have: $j\in J(\theta^{*}) \Leftrightarrow
|\hat{\theta}_{j}|>c_{2}r$. This yields
$\mathrm{sign}(\tilde{\theta}_{j})=\mathrm{sign}(\hat{\theta}_{j})=\mathrm{sign}(\theta^{*}_{j})$
if $j\in J(\theta^{*})$ on the event $\A$. If $j\not\in
J(\theta^{*})$, $\mathrm{sign}(\theta^{*}_{j})=0$ and
$\tilde{\theta}_{j}=0$ on $\A$, so that
$\mathrm{sign}(\tilde{\theta}_{j})=0$. The same reasoning holds
true simultaneously for all $\hat{\theta}\in\hat{\Theta}$ on the
event $\A$. Thus, we get the result.
\end{proof}

\section{Appendix}
We recall here some well-known results of the theory of empirical processes.
\begin{theo}[Bousquet's version of Talagrand's concentration inequality \cite{B02}]
Let $X_{i}$ be independent variables in $\X$ distributed according
to $P$, and $\F$ be a set of functions from $\X$ to $\R$ such that
$\E(f(X))=0$, $\|f\|_{\infty}\leqslant c$ and $\|f\|^{2}\leqslant
\sigma^{2}$ for any $f\in\F$. Let
$Z=\sup_{f\in\F}\sum_{i=1}^{n}f(X_{i})$. Then with probability
$1-e^{-x}$, it holds that
\begin{equation*}
Z\leqslant \E(Z)+\sqrt{2x(n\sigma^{2}+2c\E(Z))}+\frac{cx}{3}.
\end{equation*}
\end{theo}

\begin{theo}[Symmetrization theorem \cite{VW96}, p. 108]
Let $X_{1},\ldots,X_{n}$ be independent random variables with
values in $\X$, and let $\epsilon_{1},\ldots,\epsilon_{n}$ be a
Rademacher sequence independent of $X_{1},\ldots,X_{n}$. Let $\F$
ba a class of real-valued functions on $\X$. Then
\begin{equation*}
\E\left( \sup_{f\in \F}\left|
\sum_{i=1}^{n}(f(X_{i})-\E(f(X_{i}))) \right| \right) \leqslant 2
\E\left( \sup_{f\in \F} \left| \sum_{i=1}^{n}\epsilon_{i}f(X_{i})
\right|\right).
\end{equation*}
\end{theo}

\begin{theo}[Contraction theorem \cite{LT91}, p. 95]. Let
$x_{1},\ldots,x_{n}$ be nonrandom elements of $\X$, and let $\F$
be a class of real-valued functions on $\X$. Consider Lipschitz
functions $\gamma_{i}\,:\rightarrow \R$, that is,
\begin{equation*}
|\gamma_{i}(s)-\gamma_{i}(s')|\leqslant
|s-s'|,\hspace{0.25cm}\forall s,s'\in \R.
\end{equation*}
Let $\epsilon_{1},\ldots,\epsilon_{n}$ be a Rademacher sequence.
Then for any function $f^{*}\,:\,\X\rightarrow \R$, we have
\begin{equation*}
\E\left( \sup_{f\in\F}\left|
\sum_{i=1}^{n}\epsilon_{i}(\gamma_{i}(f(x_{i}))-\gamma_{i}(f^{*}(x_{i})))
\right| \right)\leqslant 2 \E\left( \sup_{f\in\F}\left|
\sum_{i=1}^{n}\epsilon_{i}((f(x_{i})-f^{*}(x_{i})) \right|
\right).
\end{equation*}
\end{theo}

\textbf{Acknowledgements:} I wish to thank Pr. Alexandre Tsybakov
for his useful advices.

\end{document}